\documentclass{article}

\usepackage[latin1]{inputenc}
\usepackage[english]{babel}
\usepackage{amsmath,amsfonts,amscd,amssymb}
\usepackage[all]{xy}
\usepackage{makeidx}
\usepackage[active]{srcltx}

\newtheorem{conj}{Conjecture}[section]
\newtheorem{thm}{Theorem}[section]
\newtheorem{cor}{Corollary}[section]
\newtheorem{lem}{Lemma}[section]
\newtheorem{prop}{Proposition}[section]
\newtheorem{deft}{Definition}[section]
\newtheorem{rem}{Remark}[section]

\newtheorem{examp}{Example}[section]

\makeindex
\title{Extending valuations to formal completions.}
\author{F. J. Herrera Govantes\footnote{Partially supported by MTM2010-19298 and FEDER.}\\
Departamento de \'Algebra\\
Facultad de Matem\'aticas\\
Avda. Reina Mercedes, s/n\\
Universidad de Sevilla\\
41012 Sevilla, Spain \\ email: jherrera@algebra.us.es \and
M. A. Olalla Acosta$^*$\\
Departamento de \'Algebra\\
Facultad de Matem\'aticas\\
Avda. Reina Mercedes, s/n\\
Universidad de Sevilla\\
41012 Sevilla, Spain\\
email: miguelolalla@algebra.us.es \and
M. Spivakovsky\\
Institut de Math\'ematiques de Toulouse\\
UMR 5219 du CNRS,\\
Universit\'e Paul Sabatier\\
118, route de Narbonne\\
31062 Toulouse cedex 9, France.\\email:
mark.spivakovsky@math.univ-toulouse.fr\and
B. Teissier\footnote{This author is grateful for the hospitality of
  the RIMS in Kyoto, where a part of this project was completed.}{}\\
Equipe `` G\'eom\'etrie et Dynamique'',\\
Institut Math\'ematique de Jussieu,\\ UMR 7586 du CNRS\\
175 Rue du Chevaleret\\
F-75013 Paris, France.\\ email: teissier@math.jussieu.fr }

\begin{document}

\newcommand{\al}{\alpha}
\newcommand{\hR}{\hat{R}}
\newcommand{\hnu}{\hat{\nu}_{-}}
\newcommand{\bnb}{\sur{B}_{\ss{n}}}
\newcommand{\End}{\textrm{End}}
\newcommand{\id}{\textrm{Id}}
\newcommand{\rg}{\mbox{rg}}
\newcommand{\rk}{\mbox{rk}}
\newcommand{\im}{\mbox{im}}
\newcommand{\ir}{ {\mathcal I} }
\newcommand{\dtr}{\mbox{degtr}}
\def\N{\mathbf{N}} 
\def\Z{\mathbf{Z}} 
\def\R{\mathbf{R}} 
\def\C{\mathbf{C}} 
\def\Q{\mathbf{Q}} 
\def\P{\mathcal{P}}
\def\F{\mathbf{F}}
\def\Tt{\mathcal{T}}
\def\sper{\mbox{Sper}}
\def\spec{\mbox{Spec}}
\def\lead{\mbox{in}}
\def\gr{\mbox{gr}}
\def\init{\mbox{in}}
\def\G{\mbox{l\hspace{-.47em}G}} 
\def\Or{\mathcal{O}}             
\def\notin{\mbox{$\in$ \hspace{-.8em}/}} 
\def\notsubset{\mbox{$\subset$ \hspace{-.9em}/}} 
\def\Ra{\Rightarrow} 
\def\Lra{\Leftrightarrow} 
\def\sgn{\mbox{sgn}} 

\def\g{{\Gamma}}
\def\h{{\Phi}}
\def\md{{\operatorname{mod}}}
\def\he{{\operatorname{ht}}}

\maketitle

\section{Introduction}

This paper is an extended version of the talk given by Miguel Olalla
at the International Conference on Valuation Theory in El Escorial in
July 2011. Its purpose is to provide an introduction to our
joint paper \cite{HOS1} without grinding through all of its technical
details. We refer the reader to \cite{HOS1} for details and proofs;
only a few proofs are given in the present paper.
\medskip

\noindent\textbf{Notation.} Throughout this paper, we will use the
following  notation:

$(R,m,k)$ a local noetherian domain

$K=QF(R)$ its field of fractions

$R_\nu$ a valuation ring dominating $R$

$\nu_{|K}\colon K^*\twoheadrightarrow\Gamma$ the restriction of $\nu$ to $K$.

$\hat{R}$ the $m-$adic completion of $R$.
\medskip

The ring $\hat R$ is local and its maximal ideal is $m\hat R$. Our
overall goal is to study extensions $\hat\nu$ of $\nu$ to $\hat R$.

The ring $\hat R$ is not in general an integral domain, so we can
only hope to extend $\nu$ to a pseudo-valuation $\hat{\nu}$ of
$\hR$. Let $P$ be the support of $\hat\nu$, that is, the prime ideal
which is mapped by $\hat\nu$ to $\infty$.

That means precisely that we want to extend $\nu$ to a valuation $\hnu$ of the
quotient $\frac{\hR}{P}$.

Such extensions $\hnu$ exists for some minimal prime ideals $P$ of $\hR$.

We shall see that $\nu$ determines a unique minimal prime $P$ of $\hR$
when $R$ is {\bf excellent}.

The purpose of our work \cite{HOS1} is to give a systematic description of all
such extensions $\hnu$, assuming that $R$ is {\bf excellent}.

Let
\begin{equation}
\Gamma\hookrightarrow\hat\Gamma\label{eq:extGamma}
\end{equation}
be an extension of ordered groups of the same rank.

Let
$(0)=\Delta_r\subsetneqq\Delta_{r-1}\subsetneqq\dots\subsetneqq\Delta_0=\Gamma$
be the isolated subgroups of $\Gamma$

and $(0)=\hat\Delta_r\subsetneqq\hat\Delta_{r-1}\subsetneqq\dots\subsetneqq
\hat\Delta_0=\hat\Gamma$
the isolated subgroups of $\hat\Gamma$, 

so that the inclusion
(\ref{eq:extGamma}) induces inclusions
\begin{eqnarray}
\Delta_\ell&\hookrightarrow&\hat\Delta_\ell\quad\text{ and}\\
\frac{\Delta_\ell}{\Delta_{\ell+1}}&\hookrightarrow&
\frac{\hat\Delta_\ell}{\hat\Delta_{\ell+1}}.
\end{eqnarray}
\medskip
Let $G\hookrightarrow\hat G$ be an extension of graded algebras
without zero divisors, such that $G$ is graded by $\Gamma_+$ and $\hat
G$ by $\hat\Gamma_+$.

\begin{deft} We say that the extension $G\hookrightarrow\hat G$ is
  \textbf{scalewise birational} if for any $x\in\hat G$ and
  $\ell\in\{1,\dots,r\}$ such that $ord\ x\in\hat\Delta_\ell$ there
  exists $y\in G$ such that $ord\ y\in\Delta_\ell$ and $xy\in G$.
\end{deft}
Of course, scalewise birational implies birational and also that
$\hat\Gamma=\Gamma$.
\medskip The
main conjecture stated in \cite{HOS1} is the following:\par
\begin{conj}\label{teissier} There exists an ideal $H$
  of $\hat R$ with $H\cap R=(0)$ and a valuation $\hat\nu_-$,
  centered at $\frac{\hat R}{H}$ and having the
  following property:\par\noindent The graded
  algebra $\gr_{\hat\nu_-}\frac{\hat R}{H}$ is a scalewise birational
  extension of $\gr_\nu R$.
\end{conj}
In \S9 we prove some partial results towards
Conjecture \ref{teissier} in the special case when $R$ is
essentially of finite type over a field. More precisely, we reduce
Conjecture \ref{teissier} to two other conjectures. Unfortunately,
the latter two conjectures remain open so far.

\section*{Acknowledgements}

We acknowledge an early precursor \cite{HeSa} of this work by
W. Heinzer and J. Sally, as well as two more recent papers \cite{CE},
\cite{CG} by S. D. Cutkosky, S. El Hitti and L. Ghezzi.

\section{Rank one valuations}

Keep the above notation. Assume that $\rk\ \nu=1$. Then the value
group $\Gamma$ is archimedian.

\medskip Let $\Phi=\nu(R\setminus\{ 0\})$. For $\beta\in\Phi$ let
$$
\begin{array}{lr} {\cal P}_\beta =&\{x\in R/\nu(x)\geq\beta\}\cr
{\cal P}^+_\beta =&\{x\in R/\nu(x)>\beta\}\end{array}
$$

We now define the main object of study of this section. Let
\begin{equation}
H:=\bigcap\limits_{\beta\in\h}(\P_\beta\hat R).\label{tag51}
\end{equation}

\begin{rem}
Since the formal completion
homomorphism $R\rightarrow\hat R$ is faithfully flat,
\begin{equation}
\P_\beta\hat R\cap R=\P_\beta\quad\text{for all }\beta\in\h.\label{tag52}
\end{equation}
Taking the intersection over all $\beta\in\h$, we obtain $H\cap R=(0)$.

In other words, we have a natural inclusion $R\hookrightarrow\frac{\hat R}H$.
\end{rem}

\begin{examp}
Let $R=k[u,v]_{(u,v)}$. Then $\hat
R=k[[u,v]]$. Consider an element 
$$w=u-\sum\limits_{i=1}^\infty c_iv^i\in\hat R\mbox{, where }c_i\in k^*,$$
transcendental over $k(u,v)$.

\medskip Consider the injective map $\varphi:k[u,v]_{(u,v)}\rightarrow
k[[t]]$ which sends $v$ to $t$ and $u$ to $\sum\limits_{i=1}^\infty c_it^i$.

\medskip Let $\nu$ be the valuation induced from the $t$-adic
valuation of $k[[t]]$ via $\varphi$.

The value
group of $\nu$ is $\mathbf Z$. For each
$\beta\in\mathbf N$,
$\P_\beta=\left(v^\beta,u-\sum\limits_{i=1}^{\beta-1}c_iv^i\right)$. 

\medskip Thus
$H=(w)$.
\end{examp}

\begin{thm}\label{th53}
\begin{enumerate}
\item The ideal $H$ is a prime ideal of $\hat R$.
\item The valuation $\nu$ extends uniquely to a valuation $\hat\nu_-$,
  centered at $\frac{\hat R}H$.
\end{enumerate}
\end{thm}

\medskip

\noindent\textit{Proof:} Let $\bar x\in\frac{\hat R}H\setminus\{0\}$. Pick a
representative $x$ of $\bar x$ in $\hat R$, so that $\bar x= x\ \md\
H$. Since $x\notin H$, we have $ x\notin \P_\alpha\hat R$ for some
$\alpha\in\h$.
\begin{lem}\label{lemma36} {\rm (See \cite{ZS}, Appendix 5, lemma
    3)} Let $\nu$ be a valuation of rank one centered in a local
  noetherian domain $(R,M,k)$. Let
$$
\h=\nu(R\setminus (0))\subset\g.
$$
Then $\h$ contains no infinite bounded sequences.
\end{lem}

\noindent\textit{Proof:} An infinite ascending sequence
$\alpha_1<\alpha_2<\dots$ in $\h$, bounded above by an element
$\beta\in\h$, would give rise to an infinite descending chain of
ideals in $\frac R{\P_\beta}$. Thus it is sufficient to prove that
$\frac R{\P_\beta}$ has finite length.

Let $\delta:=\nu(M)\equiv\min(\h\setminus\{0\})$. Since $\h$ is
archimedian, there exists $n\in\mathbf N$ such that $\beta\le n\delta$.
Then $M^n\subset \P_\beta$, so that there is a surjective map $\frac
R{M^n}\twoheadrightarrow\frac R{\P_\beta}$. Thus $\frac R{\P_\beta}$ has
finite length, as desired.
$\Box$

By Lemma \ref{lemma36}, the set $\{\beta\in\h\ |\ \beta<\alpha\}$
is finite. Hence there exists a unique $\beta\in\h$ such that
\begin{equation}
x\in \P_\beta\hat R\setminus \P_{\beta}^+\hat R.\label{tag54}
\end{equation}
Note that $\beta$ depends only on $\bar x$, but not on the choice of the
representative $ x$. Define the function $\hat\nu_-:\frac{\hat
R}H\setminus\{0\}\rightarrow\h$ by
\begin{equation}
\hat\nu_-(\bar x)=\beta.\label{tag55}
\end{equation}
By (\ref{tag52}), if $x\in R\setminus\{0\}$ then
\begin{equation}
\hat\nu_-(x)=\nu(x).
\end{equation}
It is obvious that
\begin{equation}
\hat\nu_-(x+y)\ge\min\{\hat\nu_-(x),\hat\nu_-(y)\}\label{tag57}
\end{equation}
\begin{equation}
\hat\nu_-(xy)\ge\hat\nu_-(x)+\hat\nu_-(y)\label{tag58}
\end{equation}
for all $x,y\in\frac{\hat R}H$. The point of the next lemma is to
show that $\frac{\hat R}H$ is a domain and that $\hat\nu_-$ is, in fact, a
valuation (i.e. that the inequality (\ref{tag58}) is, in fact, an equality).

\begin{lem}\label{lemma54} For any non-zero $\bar x,\bar
  y\in\frac{\hat R}H$, we have $\bar x\bar y\ne0$ and $\hat\nu_-(\bar
  x\bar y)=\hat\nu_-(\bar x)+\hat\nu_-(\bar y)$.
\end{lem}

\noindent\textit{Proof:} Let $\alpha=\hat\nu_-(\bar x)$, $\beta=\hat\nu_-(\bar
y)$. Let $ x$ and $ y$ be representatives in $\hat R$ of $\bar x$ and
$\bar y$, respectively. We have $M\P_\alpha\subset \P_{\alpha}^+$, so that
\begin{equation}
\frac{\P_\alpha}{\P_{\alpha}^+}\cong\frac{\P_\alpha}{\P_{\alpha}^++M\P_\alpha}\cong
\frac{\P_\alpha}{\P_{\alpha}^+}\otimes_Rk\cong\frac{\P_\alpha}{\P_{\alpha}^+}
\otimes_R\frac{\hat R}{M\hat R}\cong\frac{\P_\alpha\hat
R}{(\P_{\alpha}^++M\P_\alpha)\hat R}\cong\frac{\P_\alpha\hat
R}{\P_{\alpha}^+\hat R},\label{tag59}
\end{equation}
and similarly for $\beta$. By (\ref{tag59}) there exist $z\in \P_\alpha$, $w\in
\P_\beta$, such that $z\equiv x\ \md\ \P_{\alpha}^+\hat R$ and
$w\equiv y\ \md\ \P_{\beta}^+\hat R$. Then
\begin{equation}
 xy\equiv zw\ \md\ \P_{\alpha+\beta}^+\hat R.\label{tag510}
\end{equation}
Since $\nu$ is a valuation, $\nu(zw)=\alpha+\beta$, so that $zw\in
\P_{\alpha+\beta}\setminus \P_{\alpha+\beta}^+$. By (\ref{tag52}) and
(\ref{tag510}), this proves that $xy\in \P_{\alpha+\beta}\hat R\setminus
\P_{\alpha+\beta}^+\hat R$. Thus $xy\notin H$ (hence $\bar x\bar
y\ne0$ in $\frac{\hat R}H$) and $\hat\nu_-(\bar x\bar y)=\alpha+\beta$,
as desired.
$\Box$

By Lemma \ref{lemma54}, $H$ is a prime ideal of $\hat R$. By
(\ref{tag57}) and Lemma \ref{lemma54}, $\hat\nu_-$ is a valuation,
centered at $\frac{\hat R}H$. To complete the proof of Theorem
\ref{th53}, it remains to prove the uniqueness of $\hat\nu_-$. Let $x$,
$\bar x$, the element $\alpha\in\h$ and
\begin{equation}
z\in \P_\alpha\setminus \P_{\alpha}^+\label{tag511}
\end{equation}
be as in the proof of Lemma \ref{lemma54}. Then there exist
\begin{equation}
\begin{array}{rl}
u_1,\ldots , u_n&\in \P_{\alpha}^+\text{ and}\\
 v_1,\ldots ,v_n&\in\hat R
\end{array}\label{tag512}
\end{equation}
such that $ x=z+\sum\limits_{i=1}^nu_i v_i$. Letting $\bar v_i:=v_i\
\md\ H$, we obtain $\bar x=\bar z+\sum\limits_{i=1}^n\bar u_i\bar
v_i$ in $\frac{\hat R}H$. Therefore, by
(\ref{tag511})--(\ref{tag512}), for any extension of $\nu$ to a
valuation $\hat\nu '_-$, centered at $\frac{\hat R}H$, we have
\begin{equation}
\hat\nu '_-(\bar x)=\alpha=\hat\nu_-(\bar x),\label{tag513}
\end{equation}
as desired. This completes the proof of Theorem \ref{th53}.
$\Box$

\begin{deft}\label{deft55} The ideal $H$ is called the {\bf implicit prime
ideal} of $\hat R$, associated to $\nu$. When dealing with more than one
ring at a time, we will sometimes write $H(R,\nu)$ for $H$.
\end{deft}

\medskip

Recall that given a valued ring $(R,\nu)$, that is a subring $R\subseteq R_\nu$ of the valuation ring
$R_\nu$ of a valuation with value group $\Gamma$, one defines the associated graded ring
$$
\hbox{\rm gr}_\nu R=\bigoplus_{\beta\in \Gamma}\frac{{\cal P}_\beta (R)}
{{\cal P}^+_\beta (R)}=\bigoplus_{\beta\in
\Gamma_+}\frac{{\cal P}_\beta (R)}{{\cal P}^+_\beta (R)}.
$$
The second equality comes from the fact that if
$\beta\in\Gamma_-\setminus\{0\}$, we have ${\cal P}^+_\beta (R)={\cal
  P}_\beta (R)=R$.

\medskip  
  
The graded ring $G_\nu$ is the analogous one but associated to $R_\nu$ instead $R$

\begin{rem}\label{Remark56} We have the following natural
isomorphisms of graded algebras:
$$
\begin{array}{rl}
\gr_\nu R&\cong\gr_{\hat\nu_-}\frac{\hat R}H\\
G_\nu&\cong G_{\hat\nu_-}.
\end{array}
$$
\end{rem}

We will now study the behaviour of $H$ under local blowings up of $R$
with respect to $\nu$ and, more generally, under local homomorphisms.

\medskip
Let $\pi:(R,m)\rightarrow(R',m')$ be a local homomorphism of local noetherian
domains. Assume that $\nu$ extends to a rank one valuation
$\nu':K'\setminus\{0\}\rightarrow\g'$, where $\g\subset\g'$.

\medskip
The homomorphism $\pi$ induces a local homomorphism $\hat\pi:\hat
R\rightarrow\hat R'$ of formal completions. Let
$\h'=\nu'(R'\setminus\{0\})$. For $\beta\in\h'$, let $\P'_\beta$
denote the $\nu'$-ideal of $R_{\nu'}$ of value $\beta$, as above. Let
$H'=H(R',\nu')$.

\begin{lem}\label{lemma58} Let $\beta\in\h$. Then
\begin{equation}
\left(\P'_\beta\hat R'\right)\cap\hat R=\P_\beta\hat R.\label{tag516}
\end{equation}
\end{lem}

\begin{cor}\label{Corollary59} We have
$$
H'\cap\hat R=H.
$$
\end{cor}

In other words, the implicit ideals behave well under blowings up.

\begin{cor}\label{Corollary511} We have
\begin{equation}
\he\  H'\ge\he\  H.\label{tag519}
\end{equation}
In particular,
\begin{equation}
\dim\frac{\hat R'}{H'}\le\dim\frac{\hat R}H.\label{tag520}
\end{equation}
\end{cor}
It may well happen that the inequalities in (\ref{tag519}) and
(\ref{tag520}) are strict. The possibility of strict inequalities in Corollary
\ref{Corollary511} is related to the existence of subanalytic
functions, which are not analytic. We illustrate this statement by an
example in which $\he\ H<\he\ H'$.

\begin{examp} Let $k$ be a field and let
$$
\begin{array}{rl}
R&=k[x,y,z]_{(x,y,z)},\\
R'&=k[x',y',z']_{(x',y',z')},
\end{array}
$$
where $x'=x$, $y'=\frac yx$ and $z'=z$. We have $K=k(x,y,z)$, $\hat
R=k[[x,y,z]]$, $\hat R'=k[[x',y',z']]$. Let $t_1,t_2$ be auxiliary variables
and let $\sum\limits_{i=1}^\infty c_it_1^i$ (with $c_i\in k$) be an element of
$k[[t_1]]$, transcendental over $k(t_1)$. Let $\theta$ denote the valuation,
centered at $k[[t_1,t_2]]$, defined by $\theta(t_1)=1$, $\theta(t_2)=\sqrt2$
(the value group of $\theta$ is the additive subgroup of $\mathbf R$, generated
by 1 and $\sqrt2$). Let $\iota:R'\hookrightarrow k[[t_1,t_2]]$ denote the
injective map defined by $\iota(x')=t_2$, $\iota(y')=t_1$,
$\iota(z')=\sum\limits_{i=1}^\infty c_it_1^i$. Let $\nu$ denote the
restriction of $\theta$ to $K$, where we view $K$ as a subfield of
$k((t_1,t_2))$ via $\iota$. Let $\h=\nu(R\setminus\{0\})$;
$\h'=\nu(R'\setminus\{0\})$. For $\beta\in\h'$, $P'_\beta$ is generated by all
the monomials of the form ${x'}^\alpha{y'}^\gamma$ such that
$\sqrt2\alpha+\gamma\ge\beta$, together with
$z'-\sum\limits_{j=1}^ic_j{y'}^j$, where $i$ is the greatest non-negative
integer such that $i<\beta$.

Let $w':=z'-\sum\limits_{i=1}^\infty c_i{y'}^i$. Then $H'=(w')$, but
$H=H'\cap\hat R=(0)$, so that $\he\ H=0<1=\he\ H'$.
\end{examp}

\section{Introduction to the general case}

Keep the above notation. Let $r=\rk\ \nu$.

\medskip
Let
$$\g =\Delta_0\supsetneqq\cdots\supsetneqq\Delta_{r-1}\supsetneqq\Delta_r= (0)$$
be the isolated subgroups of $\g$ and
$$(0) =P_0\subsetneqq P_1\subseteqq\cdots\subseteqq P_{r-1}\subseteqq P_r= m$$
the prime $\nu$-ideals of $R$.

\medskip
For a prime ideal $P$ in $R$, $\kappa (P)$ will denote the residue
field $\frac{R_P}{PR_P}$.

Let
$$(0)=\mathbf{m}_0\subsetneqq\mathbf{m}_1\subsetneqq\dots\subsetneqq
\mathbf{m}_r=\mathbf{m}_\nu$$
be the prime ideals of the valuation ring $R_\nu$.

\medskip
By definition, our valuation
$\nu$ is a composition of $r$ rank one valuations
$\nu=\nu_1\circ\nu_2\dots\circ\nu_r$, where $\nu_\ell$ is a valuation
of the field $\kappa(\mathbf{m}_{\ell-1})$, centered at
$\frac{(R_\nu)_{\mathbf{m}_\ell}}{\mathbf{m}_{\ell -1}}$.

\subsection{Local blowings up and trees}

We consider \textit{extensions} $R\rightarrow R'$ of local rings, that
is, injective morphisms such that $R'$ is an $R$-algebra essentially
of finite type and $m'\cap R=m$. We suppose that both $R$ and $R'$ are
contained in a fixed valuation ring $R_\nu$.

Such extensions form a direct system $\{R'\}$.

\medskip
We will assume that
\begin{equation}
\lim\limits_{\overset\longrightarrow{R'}}R'=R_\nu.\label{eq:ZariskiRiemann}
\end{equation}

\begin{deft} A \textbf{tree} of $R'$-algebras is a direct system
  $\{S'\}$ of rings, indexed by the directed set
$\{R'\}$, where $S'$ is an $R'$-algebra. Note that we do not require
the maps in the direct system $\{S'\}$ to be injective. A morphism
$\{S'\}\to\{T'\}$ of trees is the datum of a map of $R'$-algebras
$S'\to T'$ for each $R'$ commuting with the tree morphisms for each
map $R'\to R''$.
\end{deft}

\begin{lem}\label{factor} Let $R\rightarrow R'$ be an extension of
  local rings. We have:\par
\noindent 1) The ideal
  $N:=m\hR\otimes_R1+1\otimes_Rm'$ is maximal in the $R$-algebra
  $\hR\otimes_RR'$.\par
\noindent 2) The natural map of completions $\hR\to\hR'$ is injective.
\end{lem}

\begin{deft} Let $\{S'\}$ be a tree of $R'$-algebras. For each $S'$,
  let $I'$ be an ideal of $S'$. We say that $\{I'\}$ is a \textbf{tree of
  ideals} if for any arrow $b_{S'S''}\colon S'\rightarrow S''$ in our
  direct system, we have $b^{-1}_{S'S''}I''=I'$. We have the obvious
  notion of inclusion of trees of ideals. In particular, we may speak
  about chains of trees of ideals.
\end{deft}

\begin{examp}
\begin{itemize}
\item The maximal ideals of the local rings of our system $\{R'\}$ form a
tree of ideals.
\item For any non-negative element $\beta\in\Gamma$, the valuation ideals
$\P'_\beta\subset R'$ of value $\beta$ form a tree of ideals of
$\{R'\}$.
\item Similarly, the $i$-th prime valuation ideals $P'_i\subset
R'$ form a tree.
\item If $rk\ \nu=r$, the prime valuation ideals $P'_i$
give rise to a chain
\begin{equation}
P'_0=(0)\subsetneqq P'_1\subseteq\dots\subseteq P'_r=m'\label{eq:treechain'}
\end{equation}
of trees of prime ideals of $\{R'\}$.
\end{itemize}
\end{examp}

Taking the limit in (\ref{eq:treechain'}), we obtain a chain
\begin{equation}
(0)=\mathbf{m}_0\subsetneqq\mathbf{m}_1\subsetneqq\dots\subsetneqq
\mathbf{m}_r=\mathbf{m}_\nu\label{eq:treechainlim}
\end{equation}
of prime ideals of the valuation ring $R_\nu$.

For each $1\leq\ell\leq r$ one has the equality
$$
\lim\limits_{\overset\longrightarrow{R'}}{\frac{R'}{P'_\ell}}=\frac{R_\nu}{\bf
  m_\ell}.
$$
Then \textbf{specifying the valuation $\nu$} is equivalent to
specifying valuations $\nu_0,\nu_1$, \dots, $\nu_r$, where $\nu_0$ is the
trivial valuation of $K$ and, for $1\le \ell\le r$, $\nu_\ell$ is a
valuation of the residue field
$k_{\nu_{\ell-1}}=\kappa(\mathbf{m}_{\ell-1})$, centered at the local
ring
$\lim\limits_{\longrightarrow}\frac{R'_{P'_\ell}}{P'_{\ell-1}R'_{P'_\ell}}=
\frac{(R_\nu)_{\mathbf{m}_\ell}}{\mathbf{m}_{\ell-1}}$ and taking its
values in the totally ordered group
$\frac{\Delta_{\ell-1}}{\Delta_\ell}$.

\medskip
From now on we fix a valuation ring $R_\nu$ dominating $R$. A local
homomorphisms $R'\rightarrow R''$ of local domains dominated by
$R_\nu$ will be called a $\nu$-\textbf{extension} of local domains.
Fix a tree $\Tt=\{ R'\}$ of noetherian local $R-$subalgebras of $R_\nu$ having
the following property: ``for each ring $R'\in\cal{T}$, all the
birational $\nu$-extensions of $R'$ belong to $\cal{T}$''.

\medskip
Moreover, we assume that $QF(R_\nu)$ equals
$\lim\limits_{\overset\longrightarrow{R'}}K'$, where $K'=QF(R')$.

\medskip
We have the following natural generalization of Conjecture \ref{teissier}:
\begin{conj}\label{teissier2} Assume that $\dim\ R'=\dim\ R$ for all
  $R'\in\mathcal{T}$.  Then there exists a tree of prime ideals $H'$
  of $\hat R'$ with $H'\cap R'=(0)$ and a valuation $\hat\nu_-$,
  centered at $\lim\limits_\to\frac{\hat R'}{H'}$ and having the
  following property:\par\noindent For any $R'\in\cal{T}$ the graded
  algebra $\gr_{\hat\nu_-}\frac{\hat R'}{H'}$ is a scalewise birational
  extension of $\gr_\nu R'$.
\end{conj}

\section{Implicit ideals}

We will define a chain of $2r+1$ prime ideals
$$
H_0\subset H_1\subset\dots\subset H_{2r}=H_{2r+1}=m\hR,
$$
satisfying $H_{2\ell}\cap R=H_{2\ell+1}\cap R=P_\ell$ for $0\le \ell\le
r$.

\medskip
We will show that these ideals $H_i$ behave well under blowings up,
that is, $H'_i\cap\hR =H_i$.

\subsection{Odd implicit ideals}

\begin{deft}
Let $0\le\ell < r$. We define our main object of study, the
$(2\ell+1)$-st implicit prime ideal $H_{2\ell+1}\subset \hR$, by
\begin{equation}
H_{2\ell+1}=\bigcap\limits_{\beta\in\Delta_\ell}
\left(\left(\lim\limits_{\overset\longrightarrow{R'}}{\cal P}'_\beta
    {\hR'}\right)\bigcap \hR\right),\label{eq:defin1}
\end{equation}
where $R'$ ranges over $\mathcal{T}$. We put $H_{2r+1}=m\hR$.
\end{deft}

We think of
(\ref{eq:defin1}) as a tree equation: if we replace $R$ by any other
$R''\in{\cal T}$ in (\ref{eq:defin1}), it defines the corresponding
ideal $H''_{2\ell+1}\subset\hat R''$.

\begin{examp}
Let $R=k[x,y,z]_{(x,y,z)}$. 

Let $\nu$ be the valuation with value group $\Gamma=\mathbf
Z^2_{lex}$, defined as follows.

Take a transcendental power series $\sum\limits_{j=1}^\infty
c_ju^j$  in a variable $u$ over $k$. 

Consider the homomorphism $R\hookrightarrow k[[u,v]]$ which sends $x$
to $v$, $y$ to $u$ and $z$ to $\sum\limits_{j=1}^\infty c_ju^j$. 

Consider the valuation $\nu$, centered at $k[[u,v]]$, defined by
$\nu(v)=(0,1)$ and $\nu(u)=(1,0)$; its restriction to $R$ will also be
denoted by $\nu$, by abuse of notation. 

We have $\nu (x)=(0,1)$, $\nu (y)= (1,0)$ and $\nu (z)=(1,0)$.

Given $\beta=(a,b)\in\mathbf Z^2_{lex}$, we have
$${\P}_\beta=x^b\left(y^a,z-c_1y-\dots-c_{a-1}y^{a-1}\right).$$ 

Then
$$\bigcap\limits_{\beta\in(0)\oplus\mathbf Z}\left({\cal P}_\beta\hat
  R\right)=(y,z)\mbox{ and }$$
$$\bigcap\limits_{\beta\in\Gamma=\Delta_0}\left({\cal P}_\beta\hat
  R\right)=\left(z-\sum\limits_{j=1}^\infty c_jy^j\right).$$ 
  
It is not
hard to show that 
$$H_1=\left(z-\sum\limits_{j=1}^\infty c_jy^j\right)\hat R\mbox{ and }$$
$$H_3=(y,z)\hat R.$$

An extension $\hat\nu$ of $\nu$ has value group $\hat\Gamma=\mathbf
Z^3_{lex}$ and is defined by $\hat\nu(x)=(0,0,1)$, $\hat\nu(y)=(0,1,0)$ and
$$\hat\nu\left(z-\sum\limits_{j=1}^\infty c_jy^j\right)=(1,0,0).$$ 
The ideal $H_1$ is the prime valuation
ideal corresponding to the isolated subgroup $(0)\oplus\mathbf
Z^2_{lex}$ of $\hat\Gamma$ and $H_3$ is the one corresponding to
$(0)\oplus(0)\oplus\mathbf Z$.

The prime ideal $H:=H_1$ and the valuation $\hnu$, induced by
$\hat\nu$ on $\frac{\hat R}H$ (that is, the valuation centered at
$\frac{\hat R}H$ with which $\hat\nu$ is composed) satisfy the
conclusion of Conjecture \ref{teissier}.

\end{examp}

\begin{examp}

Let $S=\frac{k[x,y]_{(x,y)}}{(y^2-x^2-x^3)}$. There are two distinct valuations
centered in $(x,y)$. 

Let $a_i\in k,\ i\geq 2$ be such that
$$\begin{array}{ccl}
\underbrace{\left(y+x+\sum_{i\geq 2}a_ix^i\right)} &
\underbrace{\left(y-x-\sum_{i\geq2}a_ix^i\right)} & =y^2-x^2-x^3.\\
f & g
\end{array}
$$
We shall denote by $\nu_+$ the rank one discrete valuation defined by
$$
\nu_+(x)=\nu_+(y)=1,
$$
$$
\nu_+(y+x)=2,
$$
$$
\nu_+\left(y+x+\sum_{i\geq 2}^{b-1}a_ix^i\right)=b.
$$

Now let $R=\frac{k[x,y,z]_{(x,y,z)}}{(y^2-x^2-x^3)}$. Let $\Gamma
=\mathbf Z^2_{lex}$.

Let $\nu$ be the composite valuation of the $(z)$-adic one with $\nu_+$.

The point of this example is to show that
$$
H^*_{2\ell+1}=\bigcap_{\beta\in\Delta_{\ell}}\P_{\beta}{\hat R}
$$
does not work as the definition of the $(2\ell+1)$-st implicit prime
ideal because the resulting ideal $H^*_{2\ell+1}$ is not
prime. 

Indeed, as $\P_{(a,0)}=(z^a)$, we have
$$
H_1^*=\bigcap_{(a,b)\in\Z^2}\P_{(a,b)}{\hat R}=(0).
$$
Clearly $f,g\notin H^*_1=(0)$, but $f\cdot
g=(0)$, so the ideal $H^*_1$ is not prime.

\medskip
In fact we have $H_1=(f)$ and $H_3=(z,f)$.

Let $H:=H_1$ and let $\hnu$ be the valuation of $\frac{\hat R}H\cong
k[[x,z]]$ with value group $\mathbf Z^2_{lex}$, defined by
$\hnu(x)=(0,1)$, $\hnu(z)=(1,0)$. then $H$ and $\hnu$ satisfy the
conclusion of Conjecture \ref{teissier}.

Despite what one might think from this and the previous example, there
are situations when one cannot take $H=H_1$ in Conjecture
\ref{teissier}. An explicit example of this is given in \cite{HOS1}. 

\end{examp}

\begin{prop}\label{contracts} We have $H_{2\ell+1}\cap R=P_\ell$.
\end{prop}

\begin{prop}\label{behavewell} The ideals $H'_{2\ell+1}$ behave well
  under $\nu$-extensions $R\rightarrow R'$ in $\mathcal{T}$:
  $$H_{2\ell+1}=H'_{2\ell+1}\cap \hR .$$
\end{prop}

\medskip
The main result of \cite{HOS1} is

\begin{thm}\label{odd}[Odd implicit ideals]
The implicit ideal $H_{2\ell+1}$ is prime.
\end{thm}
Next, we discuss one of the main notions used in the proof of Theorem
\ref{odd} --- that of \textbf{stable} rings in $\mathcal{T}$.

Let the notation be as above. Take an $R'\in\mathcal{T}$ and
$\overline\beta\in\frac{\Delta_\ell}{\Delta_{\ell+1}}$. Let
\begin{equation}
\P_{\overline\beta}=\left\{x\in R\ \left|\
\nu(x)\mod\Delta_{\ell+1}\ge\overline\beta\right.\right\}.\label{eq:pbetamodl}
\end{equation}
If
$$
\beta(\ell)=\min\{\gamma\in\Phi\ |\ \beta-\gamma\in\Delta_{\ell+1}\}
$$
(this makes sense because $\Phi$ is well ordered, since $R$ is
noetherian --- see \cite{ZS}, Appendix 4, Proposition 2) then
$\P_{\overline\beta}=\P_{\beta(l)}$.\par\noindent

We have the obvious inclusion of ideals
\begin{equation}
\P_{\overline\beta}\hat R\subset\P_{\overline\beta}{\hat R'}\cap
\hat R.\label{eq:stab1}
\end{equation}
A useful subtree of $\mathcal{T}$ is formed by the $\ell$-stable rings,
which we now define. An important property of stable rings, proved
in \cite{HOS1}, is that the inclusion (\ref{eq:stab1}) is an equality whenever
$R'$ is stable.
\begin{deft}\label{stable} A ring $R'\in\mathcal{T}(R)$ is said to be
  $\ell$-\textbf{stable} if the following two conditions hold:

(1) the ring
\begin{equation}
\kappa\left(P'_\ell\right)\otimes_R\left(R'\otimes_R\hat R\right)_{M'}
\label{eq:extension1}
\end{equation}
is an integral domain and

(2) there do not exist an $R''\in\mathcal{T}(R')$ and a non-trivial
algebraic extension $L$ of $\kappa(P'_\ell)$ which embeds both into
$\kappa\left(P'_\ell\right)\otimes_R\left(R'\otimes_R\hat R\right)_{M'}$
and $\kappa(P''_\ell)$.

We say that $R$ is \textbf{stable} if it is $\ell$-stable for each
$\ell\in\{0,\dots,r\}$.
\end{deft}
\begin{rem}\label{interchanging} (1) Rings of the form
  (\ref{eq:extension1}) are a basic object of study in
  \cite{HOS1}. Another way of looking at the same ring, which we often
  use, comes from interchanging the order of tensor product and
  localization. Namely, let $T'$ denote the image of the
  multiplicative system $\left(R'\otimes_R\hat R\right)\setminus M'$
  under the natural map
  $R'\otimes_R\hat R\rightarrow\kappa\left(P'_\ell\right)\otimes_R\hat
  R$. Then   the ring (\ref{eq:extension1}) equals the localization
  $(T')^{-1}\left(\kappa\left(P'_\ell\right)\otimes_R\hat R\right)$.

(2) In the special case $R'=R$ in Definition \ref{stable}, we have
$$
\kappa\left(P'_\ell\right)\otimes_R\left(R'\otimes_R\hat R\right)_{M'}=
\kappa\left(P_\ell\right)\otimes_R\hat R.
$$
If, moreover, $\frac R{P_\ell}$ is analytically irreducible then the
hypothesis that $\kappa\left(P_\ell\right)\otimes_R\hat R$ is a domain
holds automatically; in fact, this hypothesis is equivalent to
analytic irreducibility of $\frac R{P_\ell}$.

(3) Consider the special case when $R'$ is Henselian. Excellent
Henselian rings are algebraically closed inside their formal
completions, so both (1) and (2) of Definition \ref{stable} hold
automatically for this $R'$. Thus excellent Henselian local rings are
always stable.
\end{rem}
The existence of stable rings $R'\in\cal T$ is shown in \S\S7--8 of
\cite{HOS1}. More precisely, in \cite{HOS1} parallel theories
(implicit ideals, stability, etc.) are constructed not only for formal
completions, but also for henselizations and for finite local \'etale
extensions. In \S7 of \cite{HOS1} we prove the existence of stable
rings for henselization. After that, the existence of stable rings for
completion follows as an easy corollary and is proved in \S8 of \cite{HOS1}.

\medskip

The following Proposition justifies the name ``stable''.
\begin{prop}\label{largeR1} Fix an integer $\ell$, $0\le\ell\le
  r$. Assume that $R'$ is $\ell$-stable and let
  $R''\in\mathcal{T}(R')$. Then $R''$ is $\ell$-stable.
\end{prop}
The next Proposition is a technical result on which much of
\cite{HOS1} is based. For
$\overline\beta\in\frac\Gamma{\Delta_{\ell+1}}$, let
\begin{equation}
\P_{\overline\beta+}=\left\{x\in R\ \left|\
\nu(x)\mod\Delta_{\ell+1}>\overline\beta\right.\right\}.\label{eq:pbetamodl+}
\end{equation}
As usual, $\P'_{\overline\beta+}$ will stand for the analogous notion,
but with $R$ replaced by $R'$, etc.
\begin{prop}\label{largeR2} Assume that $R$ itself is $(\ell+1)$-stable and let
$R'\in\mathcal{T}(R)$.
\begin{enumerate}
\item For any $\overline\beta\in\frac{\Delta_\ell}{\Delta_{\ell+1}}$
\begin{equation}
\P'_{\overline\beta}{\hat R'}\cap\hat R=\P_{\overline\beta}\hat R.\label{eq:stab}
\end{equation}
\item For any $\overline\beta\in\frac\Gamma{\Delta_{\ell+1}}$ the natural map
    \begin{equation}\label{eq:gammaversion}
    \frac{\P_{\overline\beta}\hat R}{\P_{\overline\beta+}\hat R}\rightarrow
\frac{\P'_{\overline\beta}{\hat R'}}{\P'_{\overline\beta+}{\hat R'}}
     \end{equation}
     is injective.
\end{enumerate}
\end{prop}
\begin{cor}\label{stableimplicit} Take an integer
  $\ell\in\{0,\dots,r-1\}$ and assume that $R$ is
  $(\ell+1)$-stable. Then
\begin{equation}
H_{2\ell+1}=\bigcap\limits_{\beta\in\Delta_\ell}{\cal P}_\beta
\hat R.\label{eq:defin5}
\end{equation}
\end{cor}
\noindent\textit{Proof:} By Lemma 4 of Appendix 4 of \cite{ZS}, the ideals
$\P_{\overline\beta}$ are cofinal among the ideals $\P_\beta$ for
$\beta\in \Delta_\ell$.
$\Box$

\begin{cor}\label{stablecontracts} Assume that $R$ is stable. Take an
  element $\beta\in\Gamma$. Then $\P'_\beta{\hat R'}\cap\hat
  R=\P_\beta$.
\end{cor}
Once Proposition \ref{largeR2} and its Corollaries are established,
the proof of the primality of the odd implicit prime ideals proceeds
similarly to the case of rank one valuations, with the additional
ingredient of considering the local blowing up $R\rightarrow R'$ along
certain $\nu$-ideals.

\subsection{Even implicit ideals}

\begin{prop}\label{H2l} There exists a unique minimal prime ideal $H_{2\ell}$ of
$P_\ell \hR$, contained in $H_{2\ell+1}$.
\end{prop}
\noindent\textit{Proof:} Since $H_{2\ell+1}\cap R=P_\ell$,
$H_{2\ell+1}$ belongs to the fiber of the map $Spec\ \hat R\rightarrow
Spec\ R$ over $P_\ell$. Since $R$ was assumed to be excellent,
$S:={\hat R}\otimes_R\kappa(P_\ell)$ is a regular ring. Hence its
localization $\bar S:=S_{H_{2\ell+1}S}\cong\frac{\hat R_{H_{2\ell+1}}}{P_\ell
 \hat R_{H_{2\ell+1}}}$ is a regular {\em local} ring. In particular,
$\bar S$ is an integral domain, so $(0)$ is its unique minimal prime
ideal. The set of minimal prime ideals of $\bar S$ is in one-to-one
correspondence with the set of minimal primes of $P_\ell$, contained
in $H_{2\ell+1}$, which shows that such a minimal prime $H_{2\ell}$ is
unique, as desired.
$\Box$

\begin{cor}
$H_{2\ell}\cap R=P_\ell$.
\end{cor}

\begin{prop} We have $H_{2\ell-1}\subset H_{2\ell}$.
\end{prop}

\begin{prop}\label{behavewelleven} The ideals $H'_{2\ell}$ behave well
  under $\nu$-extensions $R\rightarrow R'$ in $\mathcal{T}$: 
  $$H_{2\ell}=H'_{2\ell}\cap \hR .$$
\end{prop}

\section{A classification of extensions of $\nu$ to $\hR$}

\begin{deft} A chain of trees $\tilde H'_0\subset\tilde
  H'_1\subset\dots\subset\tilde
H'_{2r}=m'\hR'$ of prime ideals of $\hR'$ is said to be \textbf{admissible} if:
\begin{enumerate}
\item $H'_i\subset\tilde H'_i$.
\item $\tilde H'_{2\ell}\cap R' = \tilde H'_{2\ell+1}\cap R' = P'_{\ell}$.
\item
$$
\bigcap\limits_{\overline\beta\in
\left(\frac{\Delta_\ell}{\Delta_{\ell+1}}\right)_+}
\lim\limits_{\overset\longrightarrow{R'}}
\left(\P'_{\overline\beta}\hR'+\tilde
  H'_{2\ell+1}\right)\hR'_{\tilde H'_{2\ell+2}}\cap\hR\subset\tilde
H_{2\ell+1}$$
where $\P'_{\overline\beta}$ denote the preimage in $R'$ of the
$\nu_{\ell+1}$-ideal of $\frac{R'}{P_\ell}$ of value greater than or
equal to $\overline\beta$.
\end{enumerate}
\end{deft}

\begin{thm}\label{classification} Specifying the valuation $\hnu$
  is equivalent to specifying the following data (described recursively in $i$):

(1) An admissible chain of trees $\tilde H'_0\subset\tilde
H'_1\subset\dots\subset\tilde H'_{2r}=m'\hR'$ of prime ideals of $\hR'$.

(2) For each $i$, $1\le i\le 2r$, a valuation $\hat\nu_i$ of
$k_{\hat\nu_{i-1}}$ ($\hat\nu_0$ is the trivial valuation), whose restriction to
$\lim\limits_{\overset\longrightarrow{R'}}\kappa(\tilde H'_{i-1})$ is
centered at the local ring
$\lim\limits_{\overset\longrightarrow{R'}}\frac{\hR'_{\tilde
    H'_i}}{\tilde H'_{i-1}\hR'_{\tilde H'_i}}$.

The data $\left\{\hat\nu_i\right\}_{1\le i\le 2r}$ is subject to the
following additional condition: if $i=2\ell$ is even then $rk\
\hat\nu_i=1$ and $\hat\nu_i$ is an extension of $\nu_\ell$ to
$k_{\hat\nu_{i-1}}$ (which is naturally an extension of $k_{\nu_{\ell-1}}$).
\end{thm}
\S6 of \cite{HOS1} discusses uniqueness properties of $\hnu$. Since
the statements of these results are quite technical, we omit them in
the present exposition.


\begin{thebibliography}{99}
\bibitem{CE} S. D. Cutkosky and S. El Hitti, {\em Formal prime ideals
    of infinite value and their algebraic resolution},
  Ann. Fac. Sci. Toulouse Math. (6) 19 (2010), no. 3--4, 635--649.
\bibitem{CG} S. Cutkosky and L. Ghezzi, {\em Completions of valuation
    rings}, Contemp. Math. 386 (2005), 13--34.
\bibitem{EGA} A. Grothendieck, J. Dieudonn\'e, El\'ements de
  G\'eom\'etrie alg\'ebrique, Chap. IV, Pub. Math. IHES, No.
24, 1965.
\bibitem{HeSa} W. Heinzer and J. Sally, {\em Extensions of valuations to the
completion of a local domain.} J. Pure Appl. Algebra, Vol. 71, no 2--3,
pp. 175--186, (1991).
\bibitem{HOS1} J. Herrera, M.A. Olalla, M. Spivakovsky, B. Teissier {\em
    Extending a valuation centered in a local domain to the formal
completion} Proc. London Math. Soc. (2012) 105 (3), 571-621. 
\bibitem{ZS} O. Zariski, P. Samuel {\em Commutative Algebra, Vol. II},
Springer-Verlag (1960).
\end{thebibliography}
\end{document}